\documentclass[12pt,reqno]{amsart}
	
\numberwithin{equation}{section}

\usepackage[colorlinks=true]{hyperref}
\hypersetup{urlcolor=blue, citecolor=red}

\begin{document}

\theoremstyle{plain}
\newtheorem{thm}{Theorem}[section]
\newtheorem{lem}{Lemma}[section]
\newtheorem{prop}{Proposition}[section]
\newtheorem*{cor}{Corollary}
\newtheorem*{defs}{Definitions}
\newtheorem{rem}{Remark}[section]

\newcommand{\Or}{\Omega_{r}}
\newcommand{\Orr}{\Omega_{2r}}
\newcommand{\ut}{\tilde{u}}

\newcommand{\NN}{\mathbb{N}}
\newcommand{\ZZ}{\mathbb{Z}}
\newcommand{\RR}{\mathbb{R}}
\newcommand{\CC}{\mathbb{C}}
\newcommand{\define}{\stackrel{\text{\rm def}}{=}}
\newcommand{\Ccal}{{\mathcal C}}
\newcommand{\Acal}{{\mathcal A}}
\newcommand{\Bcal}{{\mathcal B}}
\newcommand{\Vcal}{{\mathcal V}}
\newcommand{\loc}{{\text{\rm loc}}}
\newcommand{\Ecal}{{\mathcal E}}
\newcommand{\Mcal}{{\mathcal M}}
\newcommand{\Tcal}{{\mathcal T}}
\newcommand{\Ucal}{{\mathcal U}}
\newcommand{\Lcal}{{\mathcal L}}
\newcommand{\Ocal}{{\mathcal O}}
\newcommand{\rd}{\mathrm{d}}
\newcommand{\ddt}[1]{\frac{\mathrm{d}#1}{\mathrm{d}t}}
\newcommand{\dblinner}[1]{(\!({#1})\!)}

\title[Regularity of the attractor for the 2D NSE on channel-like domains]{Regularity of the global attractor for the 2D incompressible Navier-Stokes equations on channel-like domains}
\author[R. M. S. Rosa]{Ricardo M. S. Rosa}

\date{}

\address[Ricardo M. S. Rosa]{Instituto de Matem\'atica, Universidade Federal do Rio de Janeiro, Brazil}

\email[R. M. S. Rosa]{rrosa@im.ufrj.br}

\thanks{This work was supported by Coordena\c{c}\~ao de Aperfei\c{c}oamento de Pessoal de N\'{\i}vel Superior (CAPES), Brasil, grant 001, and Conselho Nacional de Desenvolvimento Cient\'{\i}fico e Tecnol\'ogico (CNPq), Bras\'{\i}lia, Brasil, grant no. 408751/2023-1.}




\begin{abstract}
The regularity of the global attractor of the incompressible Navier-Stokes equations for flows on two-dimensional domains is considered. It is assumed that domain $\Omega$ is a channel-like domain, i.e an arbitrary bounded or unbounded domain, at first without any regularity assumption on its boundary, with the only assumption that the Poincar\'e inequality holds on it. The phase space $H$ for the system is the usual closure in the $L^2(\Omega)^2$ norm of the space of smooth divergent-free vector-fields with compact support in $\Omega.$ The corresponding space obtained as the closure with respect to the $H^1(\Omega)^2$ norm is denoted by $V.$ The forcing term is assumed to belong to dual space $V'.$ It is known in this case that the global attractor exists in the phase space $H.$ It is shown in this work that the global attractor is also a compact set in $V$, and that due to the regularization effect of the equations, the solutions converge to the attractor in the norm of $V$, uniformly for initial conditions bounded in $H$. Moreover, it is shown that if the forcing term belongs to $D(A^{-s})$, for some $0<s\leq 1/2$, where $A$ is the Stokes operator, then the global attractor is compact in $D(A^{-s+1})$. If the forcing term is in $H,$ corresponding to the limit case $s=0,$ then it is further assumed that the domain is either a uniformly $\Ccal^{1,1}$ smooth domain or a bounded Lipschitz domain in order to obtain that the global attractor is compact in $D(A).$
\end{abstract}

\maketitle

\section{Introduction}

The global attractor is a fundamental mathematical object related to the permanent regimes of the flow, and many works have been devoted to studying its existence, dimension, regularity, approximation, and dynamic structure, under different sets of assumptions (see the monographs \cite{hale1988, temamdynsys, ladyzhenskaya1991, haraux1991, babinvishik1992,  sellyou2002}). In particular, the regularity of the global attractor, which is the subject of this work, is an important issue related to a detailed study of the dynamical structure of the system, to the properties of the permanent statistical regimes, and to the accuracy of numerical approximations.

The notion of global attractor for partial differential equations and, in particular, for the Navier-Stokes equations (NSE) has its origins in the work of Hopf, back in 1948, who first discussed the existence of an ``invariant manifold'' containing all the ``flows observed in the long run'' but considering instead simplified models of turbulence. Foias and Prodi \cite{foiasprodi67}, in 1967, also discussed the asymptotic behavior of two-dimensional flows in the form of invariant measures and with an embryonic form of determining modes. In 1972, Ladyzhenskaya \cite{lady72} rigorously considered the ideas proposed by Hopf and proved the existence of a compact invariant set containing all ``limiting (as $t\rightarrow \infty$) regimes.'' Then, Foias and Temam \cite{foiastemam1979}, in 1979, proved the existence of a globally attracting compact set, in the two-dimensional case, which was in fact the same invariant set constructed in \cite{lady72}, stemming from the omega-limit set of an large ball in the phase space. Foias and Temam \cite{foiastemam1979} also proved the finite-dimensionality of the attractor and the time analyticity of the solutions, with the time analyticity implying, in particular, the regularity of the global attractor as a compact set in the domain of the Stokes operator, as later observed (see e.g. \cite{fmrt2001}).

These earlier works always considered flows on bounded, smooth domains, under no-slip boundary conditions, and with the forcing term in the phase space of the system. Subsequent works have been devoted in part to extending these results to more general situations. In particular, we have proved, in \cite{rosa1998},  the existence of the global attractor in the case of a non-smooth, unbounded channel-like domain, i.e. domains in which the Poincar\'e inequality holds. Our aim in this work is to present optimal regularity results for the global attractor in this more general case of channel-like domains.

In the usual mathematical framework for the 2D Navier-Stokes equations, the phase space of the system is the space $H$ obtained as the closure, in the $L^2(\Omega)^2$ norm, of the space of smooth divergence-free vector-fields with compact support. The corresponding space obtained as the closure with respect to the $H^1(\Omega)^2$ norm is denoted by $V$, and its dual space is $V'$, with $H$ as the pivot space. 

The existence of the global attractor for the 2D Navier-Stokes equations on unbounded domains was first obtained assuming the forcing term lies in some appropriate weighted space compactly included in $H$: Abergel \cite{abergel} considered an infinitely long channel, and Babin \cite{babin} considered a channel-like domain; in both cases the boundary is assumed to be smooth, and the Poincar\'e inequality holds. Ladyzhenskaya \cite{lady92} was the first to consider the attractor on non-smooth domains, assuming however the bounded case, using different energy-type inequalities for the a priori estimates. The inequalities obtained in \cite{lady92} also allowed more general forces lying in the dual space $V'$ for the existence of the global attractor in $H$. The dimension of the attractor in the non-smooth case was addressed in \cite{bps2000}. Ilyin \cite{ilyin1996} also considered unbounded non-smooth domains, but with the restriction of having finite measure.

Then, we considered, in \cite{rosa1998}, the more general case of unbounded non-smooth channel-like domains. The forcing term was also relaxed to be in $V'$. The attractor in this case was shown to be compact in $H$ and bounded in $V$. The technique used in this more general case was different than the previous ones in that we used the energy-equation method introduced in \cite{ball, ball3} for wave-type equations on bounded domains. This was the first work adapting this technique to parabolic-like equations, which sprawled a series of works by other researchers, on different nonlinear parabolic-like equations other than the NSE.

Following our work, Ju \cite{ju2000} addressed a regularity issue, namely whether the attractor is also compact in $V$. This was proved to be true, but only in the case of smooth unbounded channel-like domains with forces in $H$. The technique used was the energy-equation method for the $V$ norm, which in this case corresponds to an enstrophy equation.

The purpose of this article is to address this regularity issue under, again, the more general assumption of a non-smooth, unbounded channel-like domain with forces in the dual space $V'$. We prove that in this case the global attractor is indeed compact in $V$ (Theorem \ref{Vcompactnessthm}), which is an optimal result, since any element $u$ in $V$ may be turned into a fixed point for an appropriate forcing term in $V',$ namely $f=\nu Au + B(u, u)$. This is a new result even in the bounded non-smooth case since Ladyzhenskaya \cite{lady92} proved only the boundedness of the global attractor in $V$.

We also study further regularity properties of the global attractor depending on the regularity of the forcing term, when the forcing term is in between the spaces $H$ and $V,$ with respect to the domain of powers of the Stokes operator. More precisely, we show that if the forcing term is in $D(A^{-s})$ for some $0 < s \leq 1/2,$ then the global attractor is compact in $D(A^{-s + 1})$ (Theorem \ref{higherregularitythm}). Again, this is an optimal result in terms of the regularity of $f.$

As for the limit case $s=0,$ corresponding to the forcing term in $H,$ it is expected that the global attractor lies in $D(A).$ This, however, seems to depend on the regularity of the domain of the Stokes operator, which depends on the regularity of the boundary. If we further assume that $\Omega$ is a smooth channel-like domain (uniformly $\Ccal^{1,1}$) or that it is a bounded Lipschitz domain, then indeed we obtain the compactness of the global attractor in $D(A)$ (Theorem \ref{DAcompactnessthm}). This is an improvement over \cite{ju2000} in the sense that they only obtained the compactness in $V,$ in the case of a smooth channel-like domain. It is also an improvement on the classical result of regularity in $D(A)$, which assumed a bounded $\Ccal^2$ smooth domain \cite{foiastemam1979}.

In this general case of non-smooth channel-like domains with forces in $V'$ or even $H,$ the energy-equation method seems to be essential for proving the \emph{existence} of the global attractor and is the only one currently available. It exploits the energy equation to prove the asymptotic compactness of the semigroup in the energy space $H$. This same method has been used for the regularity (e.g. \cite{ju2000}), using the enstrophy equation to prove compactness in the ``enstrophy'' space $V$. We want to draw attention that at least for the regularity there is a more straightforward proof. Indeed, once the attractor has been proved to be compact in some phase space, the regularity of the semigroup (as a continuous map from the phase space to a more regular space) together with the invariance of the attractor yields the compactness in a more regular space since continuous functions map compact sets into compact sets. This approach however is not necessarily easy since the proof of the regularity of the semigroup can be laborious. Note that in this approach it is not sufficient to show that each solution is more regular, it is also necessary to show that this regularization is continuous with respect to the initial conditions. This is the approach we follow, but for completeness we also include the energy-equation approach for some of the results, such as the compactness in $V,$ which in the nonsmooth case requires the approach used in \cite{lady92} instead of the enstrophy equation used in \cite{ju2000}.

We remark finally that further extensions to nonhomogeneous boundary conditions and to flows which do not vanish at infinity can also be achieved in a suitable way, by subtracting an appropriate background flow (c.f. \cite{mrw1998}).

\section{Preliminaries}
\label{prelimsec}

We consider the flow of a two-dimensional homogeneous incompressible Newtonian viscous fluid enclosed in a region $\Omega\subset\RR^2$ with boundary $\partial\Omega$ and governed by the Navier-Stokes equations. We denote by $u(x,t)\in\RR^2$ and $p(x,t)\in\RR$, respectively, the velocity and the kinematic pressure of the fluid at the point $x\in\Omega$ and at time $t\ge0$, which are then determined by the following initial-boundary value problem:
\begin{equation}
  \label{nse}
  \begin{cases}
	{\displaystyle\frac{\partial u}{\partial t}}
	 - \nu \Delta u + (u\cdot\nabla)u + \nabla p = f,
		\qquad\text{in } \Omega,
	\\
	\nabla \cdot u = 0, \qquad\text{in } \Omega,
	\\
	u=0, \qquad\text{on } \partial\Omega,
	\\
	u|_{t=0}=u_0, \qquad\text{in } \Omega,
  \end{cases}
\end{equation}
where $\nu>0$ is the kinematic viscosity of the fluid and $f=f(x)\in\RR^2$ is the external mass density of body forces, which is assumed to be time independent. 

\subsection{General settings}

The domain $\Omega$ may be an arbitrary bounded or unbounded non-empty connected open set in $\RR^2$ without any regularity assumption on its boundary $\partial\Omega$ and with the sole assumption that the Poincar\'e inequality holds on it. More precisely, we assume only the following:
\begin{equation}
  \label{prepoincare}
  \text{There exists }\lambda_1>0 \text{ such that } 
      \int_{\Omega} \varphi(x)^2 \; \mathrm{d}x
        \leq \frac{1}{\lambda_1} \int_{\Omega} |\nabla \varphi(x)|^2\; \mathrm{d}x,
        \quad \forall\, \varphi \in H_0^1(\Omega).
\end{equation}
A domain of this form is called a \emph{channel-like} domain.

The mathematical framework of \eqref{nse} is now classical \cite{lady63, lions69, temamnse}: first consider $L^2(\Omega)^2$ and $H_0^1(\Omega)^2$ endowed with the inner products
\begin{equation*}
  (u,v)=\int_{\Omega} u(x) \cdot v(x) \;\mathrm{d}x,
\end{equation*}
for $u,v \in L^2(\Omega)^2,$ and
\begin{equation*}
  \dblinner{u,v}
  = \int_{\Omega}\sum_{j=1}^2 \nabla u_{j}(x)\cdot\nabla v_{j}(x)\;\mathrm{d}x,
\end{equation*}
for $u=(u_1,u_2),v=(v_1,v_2) \in H_0^1(\Omega)^2.$ The corresponding norms are $|\cdot|=(\cdot,\cdot)^{1/2}$ and $\|\cdot\|=\dblinner{\cdot,\cdot}^{1/2}$, respectively. Thanks to \eqref{prepoincare} the norm $\|\cdot\|$ is equivalent to the usual one in $H_0^1(\Omega)^2$. Let $\Ccal_\mathrm{c}^\infty(\Omega)$ denote the space of infinitely-differentiable real-valued functions with compact support on $\Omega$. Set
\begin{align*}
  & \Vcal=\{ v \in \Ccal_\mathrm{c}^\infty(\Omega)^2; 
      \;\nabla \cdot v=0 \text{ in }\Omega\}, \\
  & H=\text{closure of } \Vcal \text{ in } L^2(\Omega)^2, \\
  & V=\text{closure of } \Vcal \text{ in } H_0^1(\Omega)^2,
\end{align*}
with $H$ and $V$ endowed with the inner product and norm of respectively $L^2(\Omega)^2$ and $H_0^1(\Omega)^2$. It follows from \eqref{prepoincare} that
\begin{equation}
  \label{poincare}
    |u|^2 \le \frac{1}{\lambda_1} \|u\|^2,\qquad \forall \,u\in V.
\end{equation}

We identify $H$ with its dual and consider the dual space $V'$ of $V$, so that
$V\subseteq H \subseteq V'$, with continuous and dense inclusions. For
simplicity, we use the notation $(\cdot,\cdot)$ also for the duality product 
between $V$ and $V'$. 

We assume that $f$ belongs to $V'$. Then, we 
consider the weak formulation of \eqref{nse}: 
\begin{equation}
  \label{weaknse}
	\ddt{}(u,v) + \nu\dblinner{u,v} + b(u,u,v) = (f,v),
		\qquad \forall \,v\in V, \;\forall\, t>0,
\end{equation}
where $b:V\times V\times V \to\RR$ is a trilinear term given by
\begin{equation}
	b(u,v,w)
	=\sum_{i,j=1}^2
		\int_{\Omega} u_{i}\frac{\partial v_j}{\partial x_i} 
                 w_{j}\; \mathrm{d}x.
\end{equation}
This term satisfies the well-known orthogonality property 
\cite{lady63, lions69, temamnse}
\begin{equation}
  \label{orthogonality}
  b(u,v,v)=0,\qquad\forall\, u,v\in V.
\end{equation}
From \eqref{orthogonality} we have also the useful skew-symmetric property
\begin{equation}
  \label{skewsymmetry}
  b(u,v,w)=-b(u,w,v),\qquad\forall\, u,v,w\in V.
\end{equation}
Important estimates for this term are based on 
the following multiplicative inequality known as
Ladyzhenskaya's inequality \cite{lady63, lady92}:
\begin{equation}
  \label{ladyineq}
  \|\varphi\|_{L^4(\Omega)} 
     \leq \frac{1}{2^{1/4}} \|\varphi\|_{L^2(\Omega)}^{1/2}
     \|\nabla \varphi\|_{L^2(\Omega)}^{1/2}, 
     \qquad \forall \,\varphi\in H_0^1(\Omega).
\end{equation}
By applying H\"older's inequality and \eqref{ladyineq} twice one obtains
\begin{equation}
  \label{estimateforb}
  |b(u,v,w)| \leq \frac{1}{2^{1/2}} |u|^{1/2} \|u\|^{1/2}\|v\|
   |w|^{1/2} \|w\|^{1/2}, \qquad \forall\, u,v,w\in V.
\end{equation}

The weak formulation \eqref{weaknse} is equivalent to the functional equation
\begin{equation}
  \label{functionalnse}
  u_t+\nu Au  + B(u)=f, \qquad \text{ in } V',
\end{equation}
where $u_t=\rd u/\mathrm{d}t$; $A:V\to V'$ is the abstract Stokes operator defined by
\begin{equation}
  (Au,v) = \dblinner{u,v}, \qquad \forall\, u,v\in V;
\end{equation}
and $B(u)=B(u,u)$ is the quadratic term associated with the bilinear operator 
$B:V\times V\to V'$ defined by
\begin{equation*}
  (B(u,v),w) = b(u,v,w), \qquad \forall\, u,v,w\in V.
\end{equation*}

The Stokes operator $A$ is an isomorphism from $V$ into $V'$. We can also
consider $A$ restricted to $D(A)=\{u\in V; \;Au\in H\}$. Such operator
$A:D(A)\subseteq H \rightarrow H$ is closed, positive and self-adjoint.
In general it does not have a compact inverse since we are not assuming
the domain to be bounded. Neither do we have the characterization of
$D(A)$ as $V\cap H^2(\Omega)^2$ since we are not assuming the boundary
$\partial\Omega$ to have any regularity. Nevertheless, we can consider
positive and negative powers $A^s$ of $A$, $s\in \RR$, and we
have at least that
\begin{equation}
  \label{powersofA}
  D(A^{1/2})=V, \quad D(A^0)=H, \quad D(A^{-1/2})=V'.
\end{equation}
We recall the embeddings
\begin{equation}
  \label{VLqemb}
  V \subseteq H_0^1(\Omega)^2 \subseteq L^p(\Omega)^2, \qquad 2\leq p < \infty.
\end{equation}
Thus, by interpolation,
\begin{equation}
  \label{fractionalemb}
  D(A^s) \subseteq L^q(\Omega)^2, \qquad 2\leq q < \frac{2}{1-2s}, 
    \quad 0 < s \leq \frac{1}{2}.
\end{equation}

From \eqref{skewsymmetry} and \eqref{estimateforb} it 
follows that $B(\cdot,\cdot)$ is continuous from $V$ into $V'$ and satisfies
\begin{equation}
  \label{estimateforB}
  \|B(u,v)\|_{V'} \leq \frac{1}{2^{1/2}} |u|^{1/2} \|u\|^{1/2}
   |v|^{1/2} \|v\|^{1/2}, \qquad \forall\, u,v\in V.
\end{equation}

Given $0<s\leq 1/2$, 
we also find using H\"older's inequality that
\[ |b(u,v,w)| \leq \|u\|_{L^p(\Omega)^2}\|v\|\|w\|_{L^q(\Omega)^2},
    \qquad \text{for} \quad 2 < q < \frac{2}{1-2s} \leq \infty, 
    \quad \frac{1}{p} + \frac{1}{q} = \frac{1}{2}.
\]
Then, using the embeddings \eqref{VLqemb} and \eqref{fractionalemb}, we obtain
\[ |b(u,v,w)| \leq c_s \|u\| \|v\| |A^s w|,
\]
for a constant $c_s$ depending on $s$. This means that
\begin{equation}
  \label{estimateforBinminuss}
  B(\cdot,\cdot):V\times V \rightarrow D(A^{-s}) \text{ is continuous for }
   0< s \leq \frac{1}{2}.
\end{equation}

We know \cite{temamnse} that given $f\in V'$ and $u_0\in H$, there exists 
a unique solution $u$ in $\Ccal([0,\infty), H)\cap L_\loc^2(0,\infty;V)$ 
of \eqref{weaknse}, hence \eqref{functionalnse} as well, with $u(0)=u_0$.
Moreover, $u_t$ belongs to $L_\loc^2(0,\infty;V')$.

Now, let $u=u(t)$, $t\geq 0$, be a solution of \eqref{weaknse}.
Since $u\in L^2(0,T;V)$ and $u_t\in L^2(0,T;V')$, we have
\begin{equation}
  \frac{1}{2} \ddt{} |u|^2 = (u_t,u),
\end{equation}
so that from \eqref{functionalnse}, and using
the orthogonality property \eqref{orthogonality}, we have that
\begin{equation}
  \label{energyeq}
  \frac{1}{2} \ddt{}|u|^2 + \nu\|u\|^2 = (f,u), \qquad\;\forall\, t>0,
\end{equation}
in the distribution sense, with $t\mapsto |u(t)|^2$
absolutely continuous on $[0,\infty)$. 

From \eqref{energyeq} and using \eqref{poincare} 
one deduces the classical estimates
\begin{equation}
  \label{l2absorbingest}
  |u(t)|^2 \leq |u_{0}|^2 e^{-\nu\lambda_1 t}
    + \frac{1}{\nu^2\lambda_1} \|f\|^2_{V'},
    \qquad \forall\, t\ge 0,
\end{equation}
and
\begin{equation}
  \label{Vestimateforu}
  \int_{0}^T \|u(t)\|^2 \, \mathrm{d}t
    \leq \frac{1}{\nu} |u_{0}|^2 + \frac{T}{\nu^2} \|f\|^2_{V'},
    \qquad \forall\, T>0.
\end{equation}

According to Ladyzhenskaya \cite{lady92}, we also have that 
for all $\varepsilon >0$, the solution $u$ belongs to 
$\Ccal([\varepsilon,\infty),V)$, and $u_t$ belongs to
$L^\infty(\varepsilon,\infty;H)\cap L_\loc^2(\varepsilon,\infty;V)$.
Moreover, $u_{tt}$ belongs to $L_\loc^2(\varepsilon,\infty;V')$, with
\begin{equation}
  \label{functionalnsesubt}
  u_{tt} + \nu A u_t + B(u,u_t) + B(u_t,u) = 0.
\end{equation}
These regularity results follow from energy-type inequalities for 
$\|u\|$ and $|u_t|$, obtained by taking (for the Galerkin approximations) 
the inner product in $H$ of $u_t$ first with
\eqref{functionalnse} and then with \eqref{functionalnsesubt},
\begin{equation}
  \label{utenergyeq}
  |u_t|^2 + \frac{\nu}{2}\ddt{} \|u\|^2 + b(u,u,u_t) = (f,u_t),
\end{equation}
and
\begin{equation}
  \label{uttenergyeq}
  \frac{1}{2}\ddt{} |u_t|^2 + \nu \|u_t\|^2 + b(u_t,u,u_t) = 0.
\end{equation}
By multiplying the first equation by $t$ and the second one by $t^2$ 
one arrives eventually at
\[ \frac{1}{2} \frac{\rd }{\rd t} \left( \nu t \|u\|^2 + t^2|u_t|^2\right)
      + \frac{\nu t^2}{2} \|u_t\|^2  
      \leq \nu \|u\|^2 + \frac{1}{\nu}\|f\|_{V'}^2
        +\frac{1}{\nu}|u|^2\|u\|^2
        +\frac{t^2}{\nu}\|u\|^2|u_t|^2.
\]
From that, the following estimates are obtained
for a solution with initial condition $u_0$ in $H$ and
forcing term $f$ in $V'$:
\begin{align}
  \label{tutestimate}
  \nu t \|u(t)\|^2 + t^2|u_t(t)|^2 
    & \leq C_1(T,\nu,\lambda_1, |u_0|,\|f\|_{V'}), 
    \qquad \forall\, T\geq t>0, \\
  \label{tsquarevxtestimate}
    \int_0^T t^2 \|u_t(t)\|^2 \;\mathrm{d}t 
       & \leq C_2(T,\nu,\lambda_1, |u_0|,\|f\|_{V'}),
    \qquad \forall\, T>0.
\end{align}
Then, it follows that
\begin{align*}
  \nu \tau \|u(t)\|^2 + \tau^2|u_t(t)|^2 
    & \leq C_1(\tau, \nu, \lambda_1, |u(t-\tau)|, \|f\|_{V'}), 
    & & \forall\; t\geq \tau > 0, \\
  \int_T^{T+T_0} \|u_t(t)\|^2 \;\mathrm{d}t & \leq 
    \frac{1}{\tau^2}C_2(T_0+\tau, \nu, \lambda_1, 
          |u(T-\tau)|, \|f\|_{V'}),
    & & \forall\, T+T_0\geq T \geq \tau>0.
\end{align*}
Thus, using \eqref{l2absorbingest} applied to $|u(t-\tau)|$ and
$|u(T-\tau)|$ we arrive at (see \cite{lady92}) 
\begin{align}
  \label{Vestimate}
  \|u(t)\| & \leq C_3(\tau, \nu, \lambda_1, |u_0|, \|f\|_{V'}), 
    \qquad \forall\, t\geq \tau>0, \\
  \label{utHestimate}
   |u_t(t)|   & \leq C_4(\tau, \nu, \lambda_1, |u_0|, \|f\|_{V'}), 
    \qquad \forall\, t\geq \tau>0, \\
  \label{utVestimate}
  \int_T^{T+T_0} \|u_t(t)\|^2 \;\mathrm{d}t & \leq 
    C_5(T_0, \tau,\nu, \lambda_1, |u_0|, \|f\|_{V'}),
    \qquad \forall\, T+T_0\geq T \geq \tau>0. 
\end{align}

We denote by $\{S(t)\}_{t\ge0}$ the semigroup in $H$ 
associated with \eqref{weaknse},
\begin{equation*}
	S(t)u_{0} = u(t), \qquad t\ge0,
\end{equation*}
where $u$ is the solution of \eqref{weaknse} with $u(0)=u_{0}\in  H$.
It follows from \eqref{l2absorbingest} that the set
\begin{equation}
  \label{Bcaldef}
  \Bcal_H=\left\{ v\in H; |v| \leq \rho_0 \equiv
    \frac{1}{\nu}\sqrt{\frac{2}{\lambda_1}} \|f\|_{V'} \right\}
\end{equation}
is absorbing in $H$ for the semigroup. It has been proved in
\cite{rosa1998} that within this framework $\{S(t)\}_{t\geq 0}$ possesses 
a compact global attractor $\Acal$ in $H$. The estimate \eqref{Vestimate} 
shows that
there exists also an absorbing set $\Bcal_V$ (in the sense that it absorbs 
the solutions uniformly for initial data bounded in $H$) which is bounded 
in $V$, say
\begin{equation}
  \label{BcalVdef}
  \Bcal_V = \left\{ v\in V; \|v\| \leq \rho_1 \equiv
     C_3(1, \nu, \lambda_1, \rho_0, \|f\|_{V'}) \right\}.
\end{equation}
As remarked in \cite{rosa1998}, the
estimate \eqref{Vestimate} together with the invariance of
$\Acal$ implies that $\Acal$ is bounded in $V$. Our aim is to
obtain further regularity properties for $\Acal$.

\subsection{More regular domains}
\label{secmoreregulardomains}

The compactness of the global attractor in $V,$ when $f\in V',$ and that in $D(A^{1 - s}),$ when $f\in D(A^{-s}),$ only assume that the Poincar\'e inequality \eqref{prepoincare} holds on the domain $\Omega.$ No regularity or boundedness is needed. But for the compactness of the global attractor in $D(A),$ when $f\in H,$ some regularity assumption seems to be needed. With this in mind, we consider two situations.

When $\Omega$ is a uniform $\Ccal^{1,1}$ bounded or unbounded two-dimensional domain (see \cite{farwigkozonosohr2009} for the precise definition), it was proved in \cite[Theorem 1.3]{farwigkozonosohr2009} that 
\begin{equation}
  \label{DAregc11}
    D(A) = V \cap H^2(\Omega)^2.
\end{equation}

By interpolation, \eqref{DAregc11} implies
\begin{equation}
  \label{DAregc11sigma}
    D(A^\sigma) = V \cap H^{2\sigma}(\Omega)^2,
\end{equation}
for every $0 \leq \sigma \leq 1$ (see the earlier work \cite{fujitamorimoto1970} for the characterization in the case of a bounded smooth domain). Then, using the Sobolev embedding $H^{2\sigma} \subset L^\infty$ for any $\sigma > 1/2,$ we have
\begin{equation}
  \label{buvwwithAsigma}
  |b(u,v,w)| \leq \|u\|_{L^\infty(\Omega)^2}\|v\|\|w\|_{L^2(\Omega)^2} \leq C|A^\sigma u| \|v\| |w|,
\end{equation}
for $u\in D(A^\sigma), \, v\in V, \, w\in H,$ for any $1/2 < \sigma \leq 1.$ This means that the bilinear operator is bounded as an operator in
\begin{equation}
  \label{BboundedAsigmaxVinH}
  B(\cdot, \cdot) : D(A^\sigma) \times V \rightarrow H,
\end{equation}
for any $1/2 < \sigma \leq 1.$

When $\Omega$ is a bounded Lipschitz two-dimensional domain, it was proved in \cite{mitreamonniaux2008} (see also the earlier work \cite{kelloggosborn1976} for convex polygons, \cite{gabeltoksdorf2022} for recent discussions in the delicate $L^p$ case, and \cite{brownshen} for the three-dimensional case) that there exists $\epsilon=\epsilon(\Omega) > 0$ such that
\begin{equation}
  \label{DApartialregsigma}
  D(A^\sigma) \subseteq \Ccal^{2\sigma - 1}(\bar\Omega), \RR^2) \subset L^\infty(\Omega)^2,
\end{equation}
for $3/4 < \sigma < 3/4 + \epsilon.$ This implies, using again \eqref{buvwwithAsigma}, that 
\begin{equation}
  \label{BboundedAsigmaxVinH2}
  B(\cdot, \cdot) : D(A^\sigma) \times V \rightarrow H,
\end{equation}
for any $3/4 < \sigma < 3/4 + \epsilon.$

There is some belief that suitable $L^\infty$ estimates similar to that used in \eqref{buvwwithAsigma} also hold in a general non-smooth domain \cite{heywood2000, heywood2014, xie1992}, but so far the attempts have not yet completely succeeded. 

\section{Descriptions of the main results}
\label{mainresultssec}

Our first aim is to prove the following result, which asserts that $\Acal$ 
is compact in $V$.
\begin{thm}
  \label{Vcompactnessthm}
  Let $\Omega\subset \RR^2$ be an open set satisfying the Poincar\'e
  inequality \eqref{prepoincare}. Assume $\nu>0$ and $f\in V'$.
  Then, the global attractor $\Acal$ in $H$ 
  of the semigroup $\{S(t)\}_{t\ge0}$ 
  associated with the evolution equation \eqref{weaknse} 
  is compact in $V$.
\end{thm}
The proof of this theorem follows directly from the compactness of
the global attractor in $H$ and the continuity of the semigroup
from $H$ into $V$ (see Section \ref{proofvcompsection}).

Following the regularity result in Theorem \ref{Vcompactnessthm}, 
we obtain further regularity of the attractor, namely that the set
$\{u_t(t); \;t\in \RR, u\in \Acal\}$ of the first time-derivative of the 
(global) solutions $u=u(\cdot)$ in the global attractor is compact in $H$. 
Due to the functional equation
\eqref{functionalnse}, this is equivalent to showing that the map 
$u\mapsto f-\nu Au-B(u)$
is continuous in $H$ on the global attractor. We state these results as 
follows:
\begin{prop}
  \label{timecompactnessprop}
  Under the assumptions of Theorem \ref{Vcompactnessthm} 
  consider the map $F(u) = f-\nu Au - B(u)$, which is defined and continuous
  from $V$ into $V'$. Then, for every $u$ in $\Acal$, 
  $F(u)$ belongs to $H$, and the restricted map 
  $F|_\Acal:\Acal\rightarrow H$ is continuous
  when $\Acal$ is endowed with the $H$ norm. 
\end{prop}

From the above result, we obtain higher regularity of the global attractor
if the forcing term $f$ is more regular. More precisely, we have the 
following result:
\begin{thm}
  \label{higherregularitythm}
  Under the assumptions of Theorem \ref{Vcompactnessthm},
  if $f$ belongs to $D(A^{-s})$ with $0< s \leq 1/2$, 
  then $\Acal$ is compact in $D(A^{-s+1})$.
\end{thm}
Note that from \eqref{powersofA} the case $s=1/2$ is precisely 
Theorem \ref{Vcompactnessthm}. The extension is for $0<s<1/2$.

The extension of Theorem \ref{higherregularitythm} to the limit case $s=0$ (and more generally to positive powers of the Stokes operator, corresponding to negative $s$) depends on the regularity of the boundary of the physical domain $\Omega$. The crucial step here is the characterization of the domain of the Stokes operator in terms of Sobolev spaces, which allows for a control of the nonlinear term by the linear, dissipative term. The boundedness of the global attractor in $D(A)$ was obtained in \cite{foiastemam1979}, but, as remarked in \cite{fmrt2001}, the same estimate and the analyticity can be used to show that the attractor is in fact compact in $D(A)$. The same argument yields the compactness of the global attractor in $D(A)$ in the case considered in \cite{ju2000}, i.e. a smooth, unbounded channel-like domain, since in this case the characterization $D(A)=V\cap H^2(\Omega)^2$ still holds (see Section \ref{secmoreregulardomains}). When $\Omega$ is a bounded Lipschitz domain, we may not have this characterization of $D(A)$ but we still have a characterization of $D(A^\gamma),$ for $3/4 < \gamma < 3/4 + \epsilon$ (see, again, Section \ref{secmoreregulardomains}), which is sufficient to bound the nonlinear term and obtain the compactness in $D(A).$ In those regards, we have the following result.

\begin{prop}
  \label{DAcompactnessprop}
  Let $\Omega\subset \RR^2$ be an domain satisfying the Poincar\'e inequality \eqref{prepoincare}.  Assume $\nu>0$ and $f\in H$. Suppose that it is possible to bound the nonlinear term so that $B:D(A^\sigma)\times D(A^\sigma) \rightarrow H$ is continuous, for some  $\sigma<1$. Then, the global attractor in $H$ of the semigroup $\{S(t)\}_{t\ge0}$ associated with the evolution equation \eqref{weaknse} is compact in $D(A)$.
\end{prop}

As a consequence of Proposition \ref{DAcompactnessprop} and the regularity results \eqref{BboundedAsigmaxVinH} and \eqref{BboundedAsigmaxVinH2} for the Stokes operator described in Section \ref{secmoreregulardomains}, we have the following corollary.
\begin{thm}
  \label{DAcompactnessthm}
  Let $\Omega\subset \RR^2$ be a domain satisfying the Poincar\'e inequality \eqref{prepoincare}. Assume $\nu>0$ and $f\in H$. Suppose further that $\Omega$ is either uniformly $\Ccal^{1,1}$ or is a bounded Lipschitz domain. Then, the global attractor in $H$ of the semigroup $\{S(t)\}_{t\ge0}$ associated with the evolution equation \eqref{weaknse} is compact in $D(A)$.
\end{thm}

Recall now that any element of a global attractor belongs to a global orbit. In the case of a bounded (not necessarily smooth) domain with $f$ in $H$ it has been proved in \cite{lady92} that for any initial condition in $V$, the corresponding solution is analytic in time on a neighborhood of the positive real axis as a function with values in $H$. This result can be easily extended (following \cite{lady92}) to the hypotheses of Theorem \ref{Vcompactnessthm}, with the solutions being analytic as functions with values in $V$. Then, the backward uniqueness property (see e.g. \cite{constantinfoias88}) holds, and hence the global orbit through a given point in the global attractor is unique. This fact is not crucial in what follows but it simplifies slightly the presentation.

\section{Proofs using the regularization of the semigroup}

We first prove the regularization properties of the semigroup, then we use them to prove the main results.

\subsection{Regularity of the semigroup}

We will show that 
for positive time $t>0$, the map $S(t)$ is continuous from $H$ into $V$.
In fact, it is a Lipschitz map. By subtracting two solution $u(t)$ and $v(t)$,
$t\geq 0$, with $u(0)=u_0$, $v(0)=v_0$, and setting $w(t)=u(t)-v(t)$, we
find the equations
\begin{align}
  \label{wteq}
  & w_t + \nu A w + B(w,u) + B(v,w) = 0, \\
  \label{wtteq}
  & w_{tt} + \nu A w_t + B(w_t,u) + B(w,u_t) + B(v_t,w) + B(v,w_t) = 0.
\end{align}
Take the inner product in $H$ of $w$ with \eqref{wteq} to find
\[ \frac{1}{2}\frac{\rd}{\rd t} |w|^2 + \nu \|w\|^2 + b(w,u,w) = 0.
\]
Using H\"older's inequality and \eqref{ladyineq} we obtain
\[ \frac{1}{2}\frac{\rd}{\rd t} |w|^2 + \nu \|w\|^2 
    \leq \frac{1}{2^{1/2}} |w|\|w\|\|u\| 
    \leq  \frac{1}{4\nu}\|u\|^2|w|^2 + \frac{\nu}{2}\|w\|^2.
\]
Thus,
\begin{equation}
  \label{wHeq}
  \frac{\rd}{\rd t} |w|^2 + \nu \|w\|^2 \leq \frac{1}{2\nu}\|u\|^2|w|^2.
\end{equation}
Dropping the second term on the left-hand side and using Gronwall's lemma
yield
\[ |w(t)|^2 \leq e^{2^{-1}\nu^{-1}\int_0^t\|u(s)\|^2\;\rd s} |w(0)|^2.
\]
Thanks to \eqref{Vestimateforu} we find
\begin{equation}
  \label{lipschitzinH}
  |u(t)-v(t)| \leq  C_6\left(T,\nu, \lambda_1, |u_0|, |v_0|, 
          \|f\|_{V'}\right)|u_0-v_0|, 
   \qquad \forall\, T\geq t \geq 0.
\end{equation}
Integrate \eqref{wHeq}, using \eqref{lipschitzinH}, to find also that
\begin{equation}
  \label{lipschitzforwsquareint}
  \int_0^T \|u(t)-v(t)\|^2 \;\rd t 
     \leq C_7\left(T,\nu, \lambda_1, |u_0|, |v_0|, \|f\|_{V'}\right)
      |u_0-v_0|^2, \qquad \forall\, T>0,
\end{equation}
for some other constant $C_7$.

Now, take the inner product in $H$ of $w_t$ with each of the equations 
\eqref{wteq} and \eqref{wtteq} to find
\begin{align*}
  & |w_t|^2 + \frac{\nu}{2}\frac{\rd}{\rd t} \|w\|^2 + b(w,u,w_t) + b(v,w,w_t) 
      = 0, \\
  & \frac{1}{2}\frac{\rd}{\rd t}|w_t|^2 + \nu \|w_t\|^2 + b(w_t,u,w_t) 
       + b(v_t,w,w_t) = 0.
\end{align*}
Multiplying these equations by $t$ and $t^2$ respectively and adding them
together yield
\begin{multline}
   \frac{1}{2} \frac{\rd }{\rd t} \left( \nu t \|w\|^2 + t^2|w_t|^2\right)
      + \nu t^2 \|w_t\|^2 + t b(w,u,w_t) + t b(v,w,w_t) \\
      + t^2 b(w_t,u,w_t) + t^2b(v_t,w,w_t) = \nu \|w\|^2. 
\end{multline}
Using \eqref{skewsymmetry} and \eqref{estimateforb} we estimate
\begin{align*}
   - t b(w,u,w_t) - t b(v,w,w_t) & \leq \frac{t}{2^{1/2}} 
    (|u|^{1/2}\|u\|^{1/2} + |v|^{1/2}\|v\|^{1/2}) |w|^{1/2}\|w\|^{1/2}
       \|w_t\| \\
    & \leq \frac{1}{2\nu}(|u|^{1/2}\|u\|^{1/2} + |v|^{1/2}\|v\|^{1/2})^2
       |w|\|w\| + \frac{\nu t^2}{4} \|w_t\|^2.
\end{align*}
Similarly,
\begin{align*}
  -t^2b(w_t,u,w_t) -t^2b(v_t,w,w_t) & \leq \frac{t^2}{2^{1/2}}
    \|u\||w_t|\|w_t\| + \frac{t^2}{2^{1/2}}
    |v_t|^{1/2}\|v_t\|^{1/2}\|w_t\||w|^{1/2}\|w\|^{1/2} \\
   & \leq \frac{t^2}{\nu}\|u\|^2|w_t|^2 
       + \frac{t}{4}|v_t|\|w\|^2 + \frac{t^4}{2\nu}\|v_t\|^2|w_t|^2
       + \frac{\nu t^2}{4} \|w_t\|^2.
\end{align*}
Putting the estimates together we find
\begin{multline}
   \frac{1}{2} \frac{\rd }{\rd t} \left( \nu t \|w\|^2 + t^2|w_t|^2\right)
      + \frac{\nu t^2}{2} \|w_t\|^2  
      \leq \left(\nu + \frac{t}{4}|v_t|\right) \|w\|^2 \\
        +\frac{1}{2\nu}(|u|^{1/2}\|u\|^{1/2} + |v|^{1/2}\|v\|^{1/2})^2 |w|\|w\|
        +\left(\frac{1}{\nu}\|u\|^2 + \frac{t^2}{2\nu}\|v_t\|^2\right)
          t^2|w_t|^2.
\end{multline}
By defining
\[ y(t) = \nu t \|w(t)\|^2 + t^2|w_t(t)|^2
\]
we find an inequality of the form
\[ \frac{\rd y(t)}{\rd t} \leq \alpha(t)\|w(t)\|^2 + \beta(t)|w(t)|\|w(t)\| 
        + \gamma(t)y(t)
\]
where, for any $T>0$, $\alpha(t)$ is bounded on $[0,T]$
thanks to \eqref{tutestimate}, $\beta(t)$ is square-integrable on $[0,T]$
thanks to \eqref{l2absorbingest} and \eqref{Vestimateforu}, 
and $\gamma(t)$ is integrable on $[0,T]$ 
thanks to \eqref{tsquarevxtestimate} and \eqref{Vestimateforu}. 
By applying the Gronwall Lemma
and using that $y(0)=0$ we find
\[ y(t) \leq \int_0^t e^{\int_s^t \gamma(\tau)\;\rd\tau}  
    \left\{\alpha(s)\|w(s)\|^2 + \beta(s)|w(s)|\|w(s)\| \right\}\;\rd s.
\]
Using \eqref{lipschitzforwsquareint} and the bounds for $\alpha(t)$
and $\beta(t)$ we obtain
\[ \nu t \|u(t)-v(t)\|^2 + t^2 |u_t(t)-v_t(t)|^2 \leq 
      C_8\left(T,\nu, \lambda_1, |u_0|, |v_0|, \|f\|_{V'}\right)|u_0-v_0|, 
   \quad \forall\, T\geq t \geq 0.
\]
for some different constant $C_8$. In particular we obtain for a yet
different constant
\begin{align}
  \label{ureglipschitz}
  \|u(t)-v(t)\| & \leq C_9\left(T,\tau, \nu, \lambda_1, |u_0|, |v_0|, 
          \|f\|_{V'} \right)|u_0-v_0|, 
     & & \forall\, T\geq t \geq \tau > 0, \\
  \label{utlipschitz}
  |u_t(t)-v_t(t)| & \leq C_{10}\left(T, \tau, \nu, \lambda_1, |u_0|, |v_0|, 
          \|f\|_{V'}\right)|u_0-v_0|, 
     & &  \forall\, T\geq t \geq \tau > 0.
\end{align}

This proves the Lipschitz continuity of $S(t)$, from bounded subsets of $H$ 
into $V$, and the Lipschitz continuity of $\rd S(t)/\rd t$, from
bounded subsets of $H$ into $H$, for all $t>0$.

\subsection{Compactness of the global attractor in \texorpdfstring{$V$}{V} - proof of Theorem \ref{Vcompactnessthm}}
\label{proofvcompsection}

We know from \cite{rosa1998} that the global attractor $\Acal$ exists in the phase space $H$, which means that it is a compact set in $H$ which is invariant by the flow. Consider then the time-one map $S(1):H\rightarrow H$.
Since the global attractor $\Acal$ is invariant, the time-one map
is a bijection when restricted to $\Acal$. On the other hand,
we have just proved that $S(1)$ is continuous from $H$ into $V$. 
Hence, $S(1)$ is a continuous
map of the set $\Acal$ endowed with the topology of $H$ onto the
set $\Acal$ endowed with the topology of $V$. Since continuous maps
take compact sets into compact sets, the set $\Acal$ endowed with the
topology of $V$ must be compact. Hence, $\Acal$ is a compact set in $V$.
This completes the proof of Theorem \ref{Vcompactnessthm}.

\subsection{Compactness in \texorpdfstring{$H$}{H} of the time-derivative of the solutions - proof of Proposition \ref{timecompactnessprop}}

Since $\Acal$ is invariant, the time-one map $S(1)$ is a bijection
in $\Acal$. Since $\Acal$ is compact in $H$, it follows that $S(1)$ is 
actually a homeomorphism of $\Acal$ in the topology of $H$. In particular, 
the inverse $S(1)|_{\Acal}^{-1}$ of $S(1)$ restricted to $\Acal$ 
exists and is continuous with respect to the topology of $H$.
Moreover, from the estimate \eqref{utlipschitz} we see that the map that takes
$u_0=u(0)$ into $u_t(1)$ is continuous from $H$ into $H$. But this
map is precisely $F(S(1)u_0)$. Finally, the map $u_0\mapsto F(u_0)$
in $\Acal$ can be written as a composition of $v_0\mapsto F(S(1)v_0)$
with $v_0=S(1)^{-1}u_0$, hence it is continuous as the composition
of two continuous functions. This completes the proof of Proposition
\ref{timecompactnessprop}.

\subsection{Higher-order regularity - proofs of Theorems \ref{higherregularitythm} and \ref{DAcompactnessthm} and of Proposition \ref{DAcompactnessprop}}

Assume the hypotheses of Theorem \ref{higherregularitythm} are valid. We write
\begin{equation}
  \label{forhigherreg}
  \nu A^{-s+1}u = A^{-s} f - A^{-s} B(u) - A^{-s} F(u),
\end{equation}
where $F(u)=f-\nu Au - B(u)$, which by Proposition \ref{timecompactnessprop} 
is continuous from $\Acal$ into $H$, 
when we consider the topology on $\Acal$ to be say $H$.

If $\{u_n\}_n$ is a sequence in $\Acal$, then, by the compactness of $\Acal$
in $V$, there exists a subsequence $\{u_{n'}\}_{n'}$ and an element
$w\in \Acal$ such that $u_{n'}\rightarrow w$ in $V$, hence also in $H$. 
By Proposition 
\ref{timecompactnessprop}, $F(u_{n'})\rightarrow F(w)$ in $H$, thus
$A^{-s} F(u_{n'})\rightarrow A^{-s} F(w)$ in $D(A^{s})$
and, in particular, in $H$. From the continuity of the bilinear term
from $V\times V$ into $D(A^{-s})$ (see \eqref{estimateforBinminuss}), 
we obtain that
$A^{-s} B(u_{n'})\rightarrow A^{-s} B(w)$ in $H$. Finally, since 
$A^{-s} f$ is assumed to belong to $H$, we obtain from 
\eqref{forhigherreg} that $A^{-s+1}u_{n'} \rightarrow A^{-s+1}w$.
This proves that $\Acal$ is bounded and compact in $D(A^{-s+1})$,
so that Theorem \ref{higherregularitythm} is proved.

In order to prove Proposition \ref{DAcompactnessprop}, we first observe
that $\Acal$ is compact in $D(A^{-s+1})$ for any $s>0$; this
follows from Theorem \ref{higherregularitythm} since $f\in H\subseteq
D(A^{-s})$ for any such $s$. Hence, in particular, $\Acal$ is compact 
in $D(A^\sigma)$, for any $\sigma<1$ allowed in the assumptions of Proposition 
\ref{DAcompactnessprop}. Consider now
\begin{equation}
  \label{forDAreg}
  \nu Au = f - B(u) - F(u).
\end{equation}
For $u\in \Acal\subseteq D(A^\sigma)$, we have that $B(u)$ 
and $F(u)$ belong to $H$. Moreover $f$ is assumed to be in $H$. Thus,
$Au\in H$, which means that $\Acal \subseteq D(A)$.
The compactness of $\Acal$ in $D(A)$ also follows from
\eqref{forDAreg}. Indeed, if $\{u_n\}_n$ is a sequence in $\Acal$,
then, since $\Acal$ is compact in $D(A^\sigma)$, it follows that
$u_{n'}\rightarrow u$ in $D(A^\sigma)$, for some $u\in \Acal$ and for
some subsequence $\{u_{n'}\}$. By the assumed continuity of $B(\cdot)$ 
from $D(A^\sigma)$ into $H$, and by the continuity of $F(\cdot)$ 
proved in Proposition \ref{timecompactnessprop}, we also have 
$B(u_{n'}) + F(u_{n'}) \rightarrow B(u) + F(u)$ in $H$. Finally, since
$f\in H$, it follows from \eqref{forDAreg} that 
$Au_{n'} \rightarrow Au$ in $H$, which proves the desired compactness.

Theorem \ref{DAcompactnessthm} just follows from Proposition \ref{DAcompactnessprop} and the continuity properties \eqref{BboundedAsigmaxVinH} and \eqref{BboundedAsigmaxVinH2}.

\section{Proofs using energy-type equations}

We now give alternative proofs of Theorem \ref{Vcompactnessthm} and of Proposition \ref{timecompactnessprop} using the energy method (see e.g. \cite{ball, ball3, ghidaglia1994, rosa1998, mrw1998, ju2000}), which yields the asymptotic compactness of the semigroup (see e.g. \cite{hale1988, temamdynsys, ladyzhenskaya1991, haraux1991, babinvishik1992,  sellyou2002}). Alternative energy-type proofs of the remaining results are likely possible, but we do not address them here.

\subsection{Compactness of the global attractor in \texorpdfstring{$V$}{V} - proof of Theorem \ref{Vcompactnessthm}) via the energy method}

Instead of using the continuity of the semigroup from $H$ into $V$
one may use the energy method to prove the compactness of $\Acal$ in
$V$. This is done as follows.

Since the global attractor $\Acal$ is an invariant set bounded
in $V$, the compactness of the global attractor in $V$ follows
as soon as we prove the \emph{asymptotic compactness}
of the semigroup $\{S(t)\}_{t\ge0}$ in the space $V$, i.e. that
$\{S(t_{n})u_{0n}\}_{n}$ is precompact in $V$ whenever
$\{u_{0n}\}_{n}$ is bounded in $V$ and $t_n\geq 0$, $t_n\to\infty$.
Indeed, if $\{v_{n}\}_n$ belongs to $\Acal$, we let $\{t_n\}$ be a
sequence of positive real numbers with $t_n\rightarrow \infty$, and,
owing to the invariance of $\Acal$, we write $v_{n} = S(t_n)u_{0n}$,
for some $u_{0n}$ in $\Acal$. Since $\Acal$ is bounded in $V$, and if
the asymptotic compactness holds in $V$, then we deduce that
$\{v_n\}_n$ is precompact in $V$. Since $\Acal$ is already
compact in $H$, it means that the limit points of $\{v_n\}_n$ belong
to $\Acal$, proving that $\Acal$ is compact in $V$.

In order to show that $\{S(t)\}_{t\ge0}$ is asymptotically
compact in $V$, we use the energy-type equations \eqref{energyeq} and
\eqref{utenergyeq}. For any solution $u=u(t)=S(t)u_0$, with $u_0\in V$,
we multiply \eqref{energyeq} by $\lambda_1$ and add it to \eqref{utenergyeq}
divided by $\nu$ to obtain, for all $t\geq 0$,
\begin{equation} 
  \label{addedenergyeq}
  \frac{1}{2}\ddt{}(\lambda_1|u|^2 + \|u\|^2) 
    + \nu\lambda_1\|u\|^2 + \frac{1}{\nu}|u_t|^2 + \frac{1}{\nu}b(u,u,u_t)
    = \lambda_1 (f,u) + \frac{1}{\nu}(f,u_t).
\end{equation}

Define
\begin{equation}
  \dblinner{u,v}_1 = \dblinner{u,v} + \lambda_1 (u,v), \qquad
  [u,v]_1=\dblinner{u,v}-\lambda_1(u,v), \qquad \forall\, u,v\in V,
\end{equation}
and set $\|u\|_1^2 = \dblinner{u,u}_1$, $[u]_1^2 = [u,u]_1$.
Clearly, $\dblinner{\cdot,\cdot}_1$, and $[\cdot,\cdot]_1$ are 
bilinear and symmetric. Moreover, from the Poincar\'e inequality 
\eqref{poincare}, we see that $\|\cdot\|_1$ is a norm in $V$ equivalent
to $\|\cdot\|$, while $[\,\cdot\,]_1$ is a seminorm in $V$.

For $u,v$ in $V$, define, also,
\[ L(u,v) = \nu\lambda_1 [u]_1^2 + \frac{2}{\nu}|v|^2,
  \qquad K(u,v) = 2\lambda_1(f,u) + \frac{2}{\nu}(f,v)
    - \frac{2}{\nu}b(u,u,v).
\]

Rewrite \eqref{addedenergyeq} as
\begin{equation}
  \label{addedmodenergyeq}
  \ddt{}\|u\|_1^2 + \nu\lambda_1 \|u\|_1^2 + L(u,u_t)
     = K(u,u_t), \qquad \forall\, t\geq 0.
\end{equation}
We now adapt the framework devised in \cite{mrw1998} to show the
asymptotic compactness in $V$ with the energy-type equation 
\eqref{addedmodenergyeq}. Let $\{u_{0n}\}_n$ be bounded in $V$ 
and $\{t_n\}_n$, $t_n \geq 0$, $t_n\to\infty$. 
By the variation of constants formula,
\begin{equation}
  \label{varconstenergyeq}
  |u_n(t)|^2 = |u_{0n}|^2 e^{-\nu\lambda_1 t}
    + \int_0^t e^{-\nu\lambda_1(t-s)}
        \left\{ K(u_n(s), u_{nt}(s)) - L(u_n(s),u_{nt}(s)) \right\} \;\rd s,
\end{equation}
where $u_n(t)=S(t)u_{0n}$, and $u_{nt}(t) = \rd u_n(t)/\mathrm{d}t$.

Since the set $\Bcal_V$ defined in \eqref{BcalVdef} is absorbing in $V$, 
there exists a time $T(B)>0$ such that
for $t_n\ge T(B)$,
\begin{equation}
  u_n(t_n) \in \Bcal_V.
\end{equation}
Thus $\{u_n(t_n)\}_n$ is weakly precompact in $V$. Since $\Acal$ is
compact in $H$ and attracts the orbits $u_n$ uniformly in $n$ as $t\rightarrow
\infty$, we find then 
\begin{equation}
  \label{Vweakconv}
  u_{n'}(t_{n'}) \rightharpoonup w\;\;
    \text{ weakly in } V \text{ and strongly in } H,
\end{equation}
for some subsequence $n'$ and for some $w\in\Acal\subset \Bcal_V$.

Similarly, for each $T\in \NN$, and passing to a further subsequence
and through a diagonalization process if necessary, we also have
\begin{equation}
  \label{VweakTconv}
  u_{n'}(t_{n'}-T) \rightharpoonup w^T\;\;
    \text{ weakly in } V \text{ and strongly in } H, \quad\forall\, T \in \NN,
\end{equation}
with $w^T\in\Acal \subset \Bcal_V$.

By the continuity of the solution semigroup $S(t)$ in $H$, we see that
$w=S(T)w^T$, for all $T\in \NN$ (the convergences of 
$\{u_{n'}(t_{n'}-T)\}_{n'}$ could now be extended 
to all $T\geq 0$, but this is not necessary).

Now, from \eqref{Vweakconv}, we find
\begin{equation}
  \|w\|_1 \leq \liminf_{n'\to\infty} \|u_{n'}(t_{n'})\|_1,
\end{equation}
and we shall now show with the help of \eqref{varconstenergyeq} that
\begin{equation*}
  \limsup_{n'\to\infty} \|u_{n'}(t_{n'})\|_1 \leq \|w\|_1.
\end{equation*}

For $T\in\NN$ and $t_n> T$ we have from \eqref{varconstenergyeq}
applied to $u_n(\cdot-T)$ that
\begin{multline}
  \label{varconstenergyeq2}
  \|u_n(t_n)\|_1^2
    = \|u_n(t_n-T)\|_1^2 e^{-\nu \lambda_1 T}
        + \int_0^T e^{-\nu \lambda_1 (T-s)}
            \left\{ K(u_n(t_n-T+s), u_{nt}(t_n-T+s)) \right. \\
            \left. - L(u_n(t_n-T+s), u_{nt}(t_n-T+s)) \right\} \;\rd s.
\end{multline}
Also, from \eqref{varconstenergyeq} applied to $w=S(T)w^T$, 
\begin{multline}
  \label{wvarconstenergyeq}
  \|w\|_1^2
    = \|w^T\|_1^2 e^{-\nu \lambda_1 T}
        + \int_0^T e^{-\nu \lambda_1 (T-s)}
            \left\{ K(S(s)w^T, \partial_s S(s)w^T) \right. \\
            \left. - L(S(s)w^T, \partial_s S(s)w^T) \right\} \;\rd s.
\end{multline}

Note that for any $T>0$,
\begin{equation}
  \label{Lfunctional}
  u \mapsto \int_0^T e^{-\nu\lambda_1(T-s)} L(u(s),u_t(s)) \;\rd s
\end{equation}
is a norm in $L^2(0,T;V)\cap W^{1,2}(0,T;H)$ equivalent to the 
usual one, hence the functional in \eqref{Lfunctional} 
is weakly lower semi-continuous in this space.
From the estimate \eqref{Vestimate} 
one finds that $\{u_{n'}(t_{n'}-T+\cdot)\}_{n'}$ 
is weakly compact in this space. Because of \eqref{VweakTconv}, we find
that this sequence converges weakly to $S(\cdot)w^T$ in this space.
Thus,
\begin{multline} 
  \int_0^T e^{-\nu\lambda_1(T-s)} L(S(s)w^T, \partial_s S(s)w^T) \;\rd s \\
     \leq \liminf_{n'\to\infty} \int_0^T e^{-\nu\lambda_1(T-s)} 
        L(u_n(t_n-T+s), u_{nt}(t_n-T+s)) \;\rd s.
\end{multline}

Now we handle the term with $K(\cdot, \cdot)$ in \eqref{varconstenergyeq2}.
From the estimates \eqref{Vestimateforu} and \eqref{utVestimate} we find 
that $\{u_{n'}(t_{n'}-T+\cdot)\}_{n'}$ converges weakly to $S(\cdot)w^T$ in
$W^{1,2}(0,T;V)$, for all $T>0$. Moreover, from the strong convergence
\eqref{VweakTconv} and the uniform continuity on $[0,T]$
of the semigroup $\{S(t)\}_{t\geq 0}$ in $H$ 
we see that $\{u_{n'}(t_{n'}-T+\cdot)\}_{n'}$ converges strongly 
to $S(\cdot)w^T$ in $L^\infty(0,T;H)$.
These are sufficient to yield
\begin{multline} 
  \int_0^T e^{-\nu\lambda_1(T-s)} K(S(s)w^T, \partial_s S(s)w^T) \;\rd s \\
     = \lim_{n'\to\infty} \int_0^T e^{-\nu\lambda_1(T-s)} 
        K(u_{n'}(t_{n'}-T+s), u_{n't}(t_{n'}-T+s)) \;\rd s.
\end{multline}

Finally, with \eqref{BcalVdef} in mind,
\begin{equation}
  \limsup_{n'\to\infty} \|u_{n'}(t_{n'}-T)\|_1^2 e^{-\nu \lambda_1 T}
        \le \rho_1^2 e^{-\nu\lambda_1 T}.
\end{equation}

Taking the $\limsup$ as $n'$ goes to infinity in \eqref{varconstenergyeq2}
and using the previous relations,
\begin{multline}
   \limsup_{n'\to\infty} \|u_{n'}(t_{n'})\|_1^2
        \leq \rho_1^2 e^{-\nu\lambda_1 T}
        + \int_0^T e^{-\nu\lambda_1(T-s)} 
         \left\{ K(S(s)w^T, \partial_s S(s)w^T) \right. \\
            \left. - L(S(s)w^T, \partial_s S(s)w^T) \right\} \;\rd s.
\end{multline}
Using \eqref{wvarconstenergyeq}, we find
\[ \limsup_{n'\to\infty} \|u_{n'}(t_{n'})\|_1^2
        \leq \rho_1^2 e^{-\nu\lambda_1 T} + \|w\|_1^2.
\]
Let $T$ go to infinity to obtain
\begin{equation}
  \label{finallimsup}
  \limsup_{n'\to\infty} \|u_{n'}(t_{n'})\|_1^2 \leq \|w\|_1^2.
\end{equation}
Since $V$ is a Hilbert space, \eqref{finallimsup} together 
with the weak convergence \eqref{Vweakconv} implies
\begin{equation}
  u_{n'}(t_{n'}) \to w\;\;
    \text{ strongly in } V.
\end{equation}
This completes the proof of Theorem \ref{Vcompactnessthm}, i.e. 
that $\Acal$ is compact in $V$.

\subsection{Compactness in \texorpdfstring{$H$}{H} of the time-derivative of the solutions - proof of Proposition \ref{timecompactnessprop}) via the energy method}

Since any element $u$ in the global attractor belongs to a global
orbit, we have $F(u)=\rd u/\mathrm{d}t$, where $\rd u/\mathrm{d}t$ is understood as
the time-derivative of this global orbit at the time it reaches $u$. 
It is then clear from
the estimate \eqref{utHestimate} that $F(u)$ belongs to $H$, 
with a uniform estimate
\begin{equation}
  \label{timederbddonatt}
  \left|\ddt{u}\right| = |F(u)| 
     \leq C_4(\tau, \nu, \lambda_1, \rho_0, \|f\|_{V'}),
  \qquad \forall \,u\in \Acal,
\end{equation}
for any given $\tau>0$.
Moreover, 
the required continuity of the map $F(\cdot)$ will follow as soon
as we prove the continuity of the time-derivatives. More precisely, suppose
that $\{w_n\}_n$ is a sequence in $\Acal$ converging in the topology 
of $V$ to some $w_0$ in $\Acal$. Since each element in $\Acal$ belongs 
to global orbit,
we can take a sequence $t_n$ of positive numbers with $t_n\rightarrow \infty$
and write $w_n=S(t_n)u_{0n}$, for some $u_{0n}$ in $\Acal$. Let $u_n=u_n(t),$
$t\in \RR,$ denote the global orbit with $u_n(0)=u_{0n}$ and $u_n(t_n)=w_n$.
Let also $w=w(t)$, $t\in \RR,$ denote the global orbit with $w(0)=w_0$.
Since $F(w_n)=\rd u_n(t_n)/\mathrm{d}t$, the continuity of $F(\cdot)$ will
follow as soon as we prove that
\[  \ddt{u_{n}(t_{n})} \rightarrow \ddt{w(0)}\;\;
    \text{ strongly in } H.
\]

Similarly to the previous section, but now using the fact that
$\Acal$ is compact in $V$ and using the backward uniqueness property
in order to simplify the calculations, we find 
that for each $T\geq 0$,
\begin{equation}
  \label{VTconv}
  u_n(t_n-T) \rightarrow w(-T)\;\;
    \text{ strongly in } V.
\end{equation}

Now, from \eqref{timederbddonatt}, we find
\begin{equation}
  \label{HweakTtimeconv}
  \ddt{u_n(t_n-T)} \rightharpoonup \ddt{w(-T)}\;\;
    \text{ weakly in } H, \qquad \forall \,T\geq 0,
\end{equation}
and, in particular,
\begin{equation}
  \left|\ddt{w(0)}\right| \leq \liminf_{n\to\infty} 
     \left|\ddt{u_n(t_n)}\right|.
\end{equation}
We shall now show with the help of the energy-type equation
\eqref{uttenergyeq} that
\begin{equation*}
  \limsup_{n\to\infty} \left|\ddt{u_n(t_n)}\right|
    \leq \left|\ddt{w(0)}\right|.
\end{equation*}

We rewrite \eqref{uttenergyeq} as
\begin{equation}
  \label{uttenergyeqmodified}
  \frac{1}{2}\ddt{} |u_t|^2 + \nu \lambda_1 |u_t|^2 + \nu [u_t]_1^2 
    + b(u_t,u,u_t) = 0.
\end{equation}
The variation of constant formula applied to this equation yields
\begin{equation}
  \label{varconsttimeenergyeq}
  |u_n'(t)|^2 = |u_{0n}'|^2 e^{-2\nu\lambda_1 t}
    - 2\int_0^t e^{-2\nu\lambda_1(t-s)}
        \left\{ b(u_n'(s), u_n(s), u_n'(s)) + \nu\lambda_1 [u_n'(s)]_1^2
        \right\} \;\rd s,
\end{equation}
where for notational simplicity we use the prime symbol $'$ for the time-derivative.

For $t_n> T>0$ we have from \eqref{varconsttimeenergyeq}
applied to $u_n(\cdot-T)$ that
\begin{multline}
  \label{varconsttimeenergyeq2}
  |u_n'(t)|^2 = |u_n'(t_n-T)|^2 e^{-2\nu\lambda_1 T} \\
    - 2\int_0^T e^{-2\nu\lambda_1(T-s)}
        \left\{ b(u_n'(t_n-T+s), u_n(t_n-T+s), u_n'(t_n-T+s))  \right. \\
   \left.       + \nu\lambda_1 [u_n'(t_n-T+s)]_1^2
        \right\} \;\rd s.
\end{multline}
Also, from \eqref{varconsttimeenergyeq} applied to $w=S(T)w^T$, 
\begin{multline}
  \label{wvarconsttimeenergyeq}
  |w'(0)|^2 = |w'(-T)|^2 e^{-2\nu\lambda_1 t} \\
    - 2\int_0^t e^{-2\nu\lambda_1(t-s)}
        \left\{ b(w'(s-T), w(s-T), w'(s-T)) 
         + \nu\lambda_1 [w'(s-T)]_1^2
        \right\} \;\rd s.
\end{multline}

Now, note that for any $T>0$,
\begin{equation}
  \label{Ltimefunctional}
  v \mapsto \int_0^T e^{-\nu\lambda_1(T-s)} [v(s)]_1^2 \;\rd s
\end{equation}
is a norm in $L^2(0,T;V)$ equivalent to the 
usual one. Hence the functional in \eqref{Ltimefunctional} 
is weakly lower semi-continuous in this space.
From the estimate \eqref{utVestimate} 
we see that $\{u_n'(t_n-T+\cdot)\}_n$ 
is weakly compact in this space. Because of \eqref{HweakTtimeconv}, we find
that this sequence converges weakly to $w'(\cdot-T)$ in this space.
Thus,
\begin{multline} 
  \int_0^T e^{-2\nu\lambda_1(T-s)} [w'(s-T)]_1^2 \;\rd s
     \leq \liminf_{n'\to\infty} \int_0^T e^{-2\nu\lambda_1(T-s)} 
      [u_n'(t_n-T+s)]_1^2 \;\rd s.
\end{multline}

Now we handle the nonlinear term in \eqref{varconsttimeenergyeq2}.
We want to prove that 
\begin{multline} 
  \label{nonlineartermconvergence}
  \int_0^T e^{-2\nu\lambda_1(T-s)} 
    b(u_n'(t_n-T+s), u_n(t_n-T+s), u_n'(t_n-T+s)) \;\rd s \\
     = \lim_{n'\to\infty} \int_0^T e^{-2\nu\lambda_1(T-s)} 
       b(w'(s-T), w(s-T), w'(s-T)) \;\rd s.
\end{multline}
We know already, from the compactness of $\Acal$ in $H$, that
$u_n(t_n-T+\cot)$ converges uniformly in $V$ to $w(\cdot)$ on bounded
time intervals. We also have from the estimates \eqref{utHestimate} and 
\eqref{utVestimate} that
that $\{u_n'(t_n-T+\cdot)\}_n$ is weakly compact in $L^2(0,T;V)$
and weak-star compact in $L^\infty(0,T;H)$. Moreover, from
the equation \eqref{functionalnsesubt} and using also the 
inequality \eqref{estimateforB} and the estimate \eqref{l2absorbingest},
we find that $\{u_n''(t_n-T+\cdot)\}_n$ is weakly compact in $L^2(0,T;V')$. 
By the Lions-Aubin Compactness Theorem 
it follows that $\{\tau_r u_n'(t_n-T+\cdot)\}_n$ is compact
in $L^2(0,T;H)$ for any truncation $\tau_r$ to compact support in a
bounded set of radius $r$ from the origin. From the weak convergence 
\eqref{HweakTtimeconv}. By triangulation and by approximating
$w(t_n-T+\cdot)$ by a compactly supported function, one may
proceed as similarly done in the proof of existence of weak solutions, 
to obtain the desired convergence \ref{nonlineartermconvergence}.

Taking the $\limsup$ as $n$ goes to infinity in \eqref{varconsttimeenergyeq2}
and using the previous relations and the estimate \eqref{timederbddonatt},
we obtain
\begin{multline}
   \limsup_{n'\to\infty} |u_n'(t_n)|^2
        \leq C_4^2 e^{-2\nu\lambda_1 T} \\
        - 2\int_0^T e^{-\nu\lambda_1(T-s)} 
         \left\{ b(w'(s-T), w(s-T), w'(s-T)) + \nu\lambda_1 [w'(s-T)]_1^2
        \right\} \;\rd s.
\end{multline}
Using \eqref{wvarconsttimeenergyeq}, we find
\[ \limsup_{n\to\infty} |u_n'(t_n)|^2
        \leq C_4^2 e^{-2\nu\lambda_1 T} + |w'(0)|^2.
\]
Let $T$ go to infinity to obtain
\begin{equation}
  \label{finaltimelimsup}
  \limsup_{n\to\infty} |u_n'(t_n)|^2 \leq |w'(0)|^2.
\end{equation}
Since $H$ is a Hilbert space, \eqref{finaltimelimsup} together 
with the weak convergence \eqref{HweakTtimeconv} implies
\begin{equation}
  u_n'(t_n) \rightarrow w'(0)\;\;
    \text{ strongly in } H.
\end{equation}
This completes the proof of Proposition \ref{timecompactnessprop}.

\end{document}